\documentclass[12pt]{article}
\usepackage{amsmath}
\usepackage{amssymb}
\usepackage{latexsym}
\usepackage{tikz}
\usepackage{hyperref}
\usepackage{epsfig}
\usepackage{graphicx}

\begin{document}
\begin{center}
%Citation:A. Hajji and D. Ziou.  Pythagorean Means for Data Selection. Technical report, Dept. informatique, Université de Sherbrooke, Qc, Canada, 2021.\\
%\line(1,0){400}\\
%[2mm]
\begin{Large}
\textbf{Pythagorean Centrality for Data Selection}\\
\end{Large}
Djemel Ziou \\
%A. Hajji and D. Ziou\\
%\line(1,0){400}\\
%[5mm]
Département d'informatique\\
Université de Sherbrooke\\
Sherbrooke, QC, Canada J1K 2R1\\
Djemel.Ziou@usherbrooke.ca
\vspace{5mm}
\end{center}
\begin{center}
\textbf{Abstract}
\end{center}
\vspace{-2mm}
 This paper provides an overview of the Pythagorean centrality  measures, which  are the arithmetic, geometric, and harmonic means. Both the evolution of their meaning through history and their geometrical interpretation are outlined. Relevant examples of use cases for each of them are introduced, spanning a variety of areas of knowledge. Their differences and similarities  are explored. Finally, the issue of which mean to use in different situations in order to make advantageous predictions is addressed.  

\section{Introduction}

The concept of the mean of a set of measurements was known in civilizations before our era. Its meaning has evolved over time,  according to needs. In ancient India, Rituparna estimated the number of leaves on a tree based on a single twig, which he multiplied by the estimated number of twigs on the branches. This can  be seen as an intuitive predecessor to the arithmetic mean, as the leaf count of a single twig was taken to be representative of any other twig~\cite{PAPPUS}. During the  civilization of the ancient Greeks, the Pythagoreans found several mathematical formulations of the mean which can give different values for the same two numbers~\cite{EISENHART}. The 9th-century mathematician  Al-Kindi introduced  cryptography  which used  several statistical concepts including  the relative frequency and the arithmetic mean~\cite{BROEMELING}. The 16th and 17th centuries brought  a shift in paradigm that led to the arithmetic mean being viewed as the "true" mean in the eyes of most people~\cite{EISENHART}. 
 
In our current knowledge,  given the non-negative measurements $x_1 < \cdots < x_n$, 
the mean can be any value between $x_1$ and $x_n$~\cite{BORWEIN}. However for practical reasons, a requirement is mandatory in order to target one value among all the values. The requirement is often expressed by informal arguments, such as the mean  value must greatly resemble  other measurements, or that it must be the closest to the others, or the one that expresses a trend of the phenomenon, or the value that compensates for errors such as  the additive noise in a signal.
No matter how the informal requirement is expressed,  the mean value  is the  closest to  the measurements, up to some transform,  according to some criterion, e.g. the value $\mu$ minimizing the quadratic error $\sum_{i=1}^{n} (T(x_i)-O(\mu))^2$, where $T()$ and $O()$ are some transforms. One may think that in order to sidestep  an arbitrary choice of a value as a mean among all, it is enough to   determine the value that satisfies a predefined criterion. But it is not so, because there can be an infinity of criteria.  In other words, we are not escaping arbitrariness altogether, but moving from a type of arbitrariness  where the mean is in $[x_1, x_n]$  to  an agreed arbitrariness where the mean value fulfils some criterion chosen arbitrarily among many others. Despite the arbitrariness, 
the mean is a concept that many among us claim to know and is the basis of reasoning in various fields. 
It is widely used for purposes such as determining  the price evolution of consumer products,  estimating  the Human Development Index, evaluating  students learning, measuring  the daily air temperature in a city, and calculating the capacitance of several capacitors connected in series. The goal of this paper is to examine  Pythagorean centrality measures: the arithmetic mean (AM), the geometric mean (GM), and the harmonic mean (HM). For this purpose, we will describe these  three means as a finite set of non-negative measurements and explore the relationship between them,  their history, and their use. Moreover, we will propose a  geometric interpretation of the three means in terms of dimensions of hyperrectangles  as well as  a statistical interpretation in terms of data selection. The next section is devoted to the definition and the usefulness  of Pythagorean centrality measures. Their history and the comparison between them are described in Section 3.  In Section 4, we present the means as a mechanism of data selection. Examples of  use cases are in section 5.

\section{Pythagorean centrality measures}
%\subsection{Arithmetic mean}
The term Pythagorean centrality measures encompass the arithmetic, geometric, and harmonic means. In all this study, we limit ourselves to the case of non-negative measurements.
Given $n$ scalar measurements $ x_1  < \cdots < x_n$, the interpretation  given our current knowledge would suggest that a Pythagorean centrality measure is located between $x_1$ and $x_n$~\cite{BORWEIN} and  is related to compensation of errors, tendency, balance, and representativeness.  These qualifications can be considered synonymous because all  can be couched as the value $\mu$ minimizing a criterion based on the squared error; i.e. $\sum_{i=1}^{n} w_i (T(x_i)-O(\mu))^2$, where $T()$ and $O()$ are some transforms. The non-negative weight $w_i$ of $x_i$ can be interpreted as its relative frequency, the empirical probability of observing it, and its relevance according to a certain criterion.   In what follows, we will use centrality, trend, and mean interchangeably and focus only on the mean of a finite set of measurements. In reality,  the geometric (resp.  harmonic) mean is also an  arithmetic mean of measurements transformed by logarithmic (resp. inverse) transform as  shown in Table~\ref{Pythagore}.    In both these cases,  the measurements in the original space are transformed, the arithmetic mean is calculated in the transformed space, and then the mean is  transported back to the original space.  The squared-error-based criterion to be minimized is expressed in the transformed space. As we will see later, this formalization of the Pythagorean centrality measures in Table~\ref{Pythagore} makes it possible to explain one of the data selection mechanisms.    
\begin{table}
\begin{center}
\begin{tabular}{|l||l|l|l|} \hline
 & AM  & GM & HM  \\ \hline
   Criterion &  \multicolumn{3}{c|}{$min_a \sum_{i=1}^{n} w_i (y_i-a)^2$, $w_i\ge0$ $\forall i$, and $\sum_{i=1}^{n} w_i=1$}   \\ \hline
Data transform & y=x & $y=ln(x)$  & $y=1/x$   \\ \hline
 Parametrization & $a =\mu$ & $a=ln(\mu)$  & $a=1/\mu$   \\ \hline
 Mean ($\mu$) & $\sum_{i=1}^{n} w_i x_i$ & $\prod_{i=1}^{n} x_i^{w_i}$ & $1/(\sum_{i=1}^{n} w_i/x_i)$   \\ \hline
\end{tabular}
\end{center}
\label{Pythagore}
\caption{Pythagorean centrality measures: the minimization problem by which each is derived, the data transform, the parametrization used, and the mean.}
\end{table}
%Another lesson to be learned from Table 1 is that the interpretation of Pythagorean  means using data transformations suggests that there are an infinity of possible transformations and therefore an infinity of means, some of which might be of practical interest. This issue will be dealt with later. 

%In ancient Greece, the arithmetic mean was used to estimate the number of crews per ship, among many other uses. Even, there are also traces of its use in ancient India for the estimation of tree leaves. The 

Most humans  default to the arithmetic mean  when asked to find the average of a sample. During an undergraduate  class  in computer science,   30 students were asked to calculate the mean of three numbers. As expected, all calculated the arithmetic mean and no one thought of using another measure of centrality.  There are countless examples of everyday use for the arithmetic mean, such as the estimation of  marks obtained by students,  the daily temperature,  the  position  of a centre of mass~\cite{TAYLOR}, and the Consumer Product Index (CPI)~\cite{STATCAN: CHAPTER 6}. It is easy to compute and  understand,   and,  unlike the  other two means, it does not require a transformation of  measurements.  The arithmetic mean is suitable for indicating the tendency when the dispersion of measurements is not high.  For the compensation of errors, a noisy measurement can be replaced by the arithmetic mean estimated from the closest measurements in the case of uncorrelated and additive noise~\cite{Ziou98}. However, it is sensitive to extreme values. For example, let's assume that the five  equally-weighted measurements of the daily temperature  reported are 8, 13, 14, 10, and 1000 degrees. Although 1000 degrees is clearly a mistake, it influences the average daily temperature, resulting in an  unrealistic value: 209 degrees.  This sensitivity of the arithmetic mean   may explain why the geometric mean has regained popularity in  many areas of knowledge. In physics, it expresses  the relationship between the classical and the relativistic Doppler effect~\cite{BAIRD}. In image processing, a geometric filter based on the geometric mean  is used for the reduction of  speckle in radar images~\cite{GEOMETRIC MEAN: SAR}. As a filter, the geometric mean  weakens the contribution of high values heavily  corrupted by multiplicative noise.  In finance, the geometric mean is used to characterize the growth rates of an investment~\cite{HULL}. It is also used to assess the effectiveness of a vaccine~\cite{Earle}.  The CPI is estimated  using the geometric mean for certain products, as  will be discussed below.  The harmonic
mean is used for  change detection in radar images~\cite{HARMONIC FILTER}, estimation of the average velocity over a trip~\cite{Falk},    averaging  the financial multiples with price in the numerator~\cite{Matthews06}, estimation of the pollutant load that can be permitted in a water quality limited stream~\cite{Limbrunner00}, and estimation of the properties of a heterogeneous system of porous media~\cite{Limbrunner00}.   
In the case of multiples, the harmonic mean can provide more accurate estimates than the arithmetic mean. This surprising result has been addressed several times in litigation in the United States~\cite{Matthews21}. 
The geometric and harmonic means favour small values and are therefore less sensitive than the arithmetic mean to  larger outliers.  In the example presented above the 
arithmetic mean is 209 degrees while the  geometric and harmonic means are around 27.02 and  13.36 degrees, respectively.  While the former is far from the median, which is 13 degrees, the latter is close to it.

\begin{figure}[]
        \centering
        \includegraphics[scale=0.3]{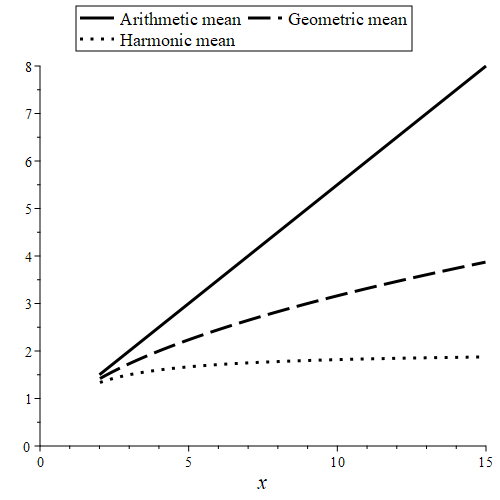}
        \vspace{-5mm}
        \caption{The three means of 1 and x.}
        \label{fig2}
    \end{figure}
    
Finally, the three means can be ranked. To illustrate this, the three means of 1 and $x$ are plotted in Fig. \ref{fig2} as a function of $x \in [2,15]$. It can be observed that  the means are contained in the interval $[0,x]$ and  AM $>$ GM $>$ HM. The inequality is valid for any finite number of measurements.  The reader can find the formal proof in~\cite{AOPS}.

\section{Evolution of the Pythagorean centrality measures}
 Note that historical clarification is required  to summarize what was defined by Pythagoras,  his disciples and what has been introduced afterwards. For example, in~\cite{VOGEL}, it is stated that the geometric mean was introduced at the beginning of the 19th century. For this purpose, the history of the three means and the current geometrical interpretation are the focus of this section.
\subsection{History}
 The arithmetic, geometric, and harmonic means of a sample appeared   during the last centuries BCE. The Babylonians used the arithmetic mean in astronomy~\cite{BABYLON},  classical Greek music theorists employs the geometric mean in tuning instruments~\cite{Huffman}, and  the harmonic mean is the basis of a solution of a division of loaves of bread by 10~\cite{Spalinger}. However, it was in ancient Greece that the study of the mean figured prominently because Greek mathematics was primarily concerned with proportions, as  is the case with ratios and the Pythagorean theorem.  The remarkable contributions of Pythagoras in the study of proportions probably led to  the three measures of centrality attributed  to Pythagoras and his disciples, i.e. Pythagorean means~\cite{Faradj04}. Much later, Al-Kindi introduced cryptography in the 9th century, using several statistical concepts, including relative frequency and arithmetic mean~\cite{BROEMELING}. In the 11th century, Al-Biruni used the midrange as a  mean  in astronomy~\cite{EISENHART}. Until the end of the Middle Ages, the notion of mean had several interpretations. One of these was the midrange~\cite{MIDRANGE}, which  is the value located equidistant from the two extreme values. This is not easily generalizable to the case of more than two measurements, which is why the definition of the arithmetic mean of a sample as we know it today was not proposed until the 16th century~\cite{EISENHART}. There is no consensus on the origin of the arithmetic mean of a population (i.e. the mathematical expectation). Some date it to the 17th century, attributing it to Blaise Pascal. Regardless, it was only formalized in the 19th century,  by Laplace~\cite{Laplace}.  Also in the Middle Ages, the harmonic mean was used to solve what is commonly called the cistern mathematical problem that appeared in Europe and India~\cite{PATBLOG}. Details regarding the development of the harmonic mean after this time period seem to be lacking in the literature. The geometric mean as we know it was not proposed until the beginning of the 19th century~\cite{VOGEL}. Surprisingly, the population geometric mean was only introduced in 2013~\cite{VOGEL}. Note that the geometric and arithmetic means have been combined to form the arithmetic-geometric mean used in the calculation of  integrals~\cite{GEO-ARI,Abramowitz}. As its name suggests, this average consists of iteratively calculating the arithmetic and geometric mean of the measurements.
 
% \subsection{AM-GM-HM inequality}

\subsection{Geometric interpretation}
\begin{figure}[]
        \centering
        \includegraphics[scale=0.6]{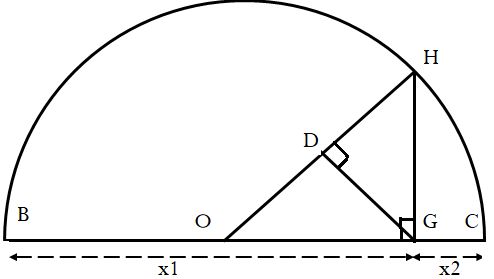}
        \caption{Representation of the three means of two numbers.}
        \label{fig1}
\end{figure}

The geometric interpretation is only provided for equal weights measurements.
Let us consider a sample of only two values, $x_1$ and $x_2$ with equal weights. The three means can be represented geometrically using a circle of diameter $x_1 + x_2$, as depicted in Fig.~\ref{fig1}. 
The lengths of BG and GC are  $x_1$ and $x_2$, respectively. The arithmetic mean coincides with the radius of the circle (OH),   HG corresponds to the geometric mean, and HD to the harmonic mean~\cite{Alsina06}. This relationship among the three means is only applicable to two numbers. We suggest an alternative geometric interpretation that can be applied to more than two numbers $( n \geq 2)$.  Let us consider the hyperrectangle $R_n$, where $x_1$, ..., $x_n$ are its dimensions. In $R_n$, there is a total of $2^{(n-1)}n$ edges. As a result, each edge is repeated $2^{(n-1)}$ times. For instance, in $R_3$, each one of the 3 edges would be repeated 4 times. The hyperperimeter of $R_n$ is given by:
\begin{equation}
    P_n = 2^{(n-1)}(x_1 + x_2 +...+x_n)
\end{equation}
According to this formula, the hyperperimeter is the arithmetic mean multiplied by a factor whose value is contingent on the number of dimensions: 
\begin{equation}
    AM = \frac{P_n}{2^{(n-1)}n}     
\end{equation}
The hypervolume $V_n = x_1 x_2 ... x_n$ of $R_n$ is the geometric mean raised to the n\textsuperscript{th} power:
\begin{equation}
    GM = \sqrt[n]{V_n}     
\end{equation}
The hyperrectangle  $R_n$ is made up of $n$ pairs of equal $(n-1)$-dimensional hyperrectangles. For example, in $R_3$, there are 6 rectangles that are equal to one another two by two. Let us consider the $j^{th}$ hyperrectangle $R_{n-1,j}$ with hypervolume  $V_{n-1,j}=\prod\limits_{i \neq j}^{n} x_i$. The arithmetic mean of the hypervolumes of the $2n$  hyperrectangles $R_{n-1,j}$ is:
\begin{equation}
    A_n =  \frac{1}{2 n} \sum^{2n}_{j = 1} V_{n-1,j} = \frac{1}{n} \sum^{n}_{j = 1} V_{n-1,j}
\end{equation} 
The harmonic mean is expressed in terms of a ratio between the hypervolume of $R_n$ and  the arithmetic mean of the hypervolumes of these n $R_{n-1}$: 
\begin{equation}
    H_n = \frac{V_n}{A_n}
\end{equation} 
In other terms, this is the ratio between the 
geometric mean of edges to the power of $n$  and the arithmetic mean $A_n$ of the hypervolumes of $n$ hyperrectangles of dimension $n-1$. To conclude, the geometric interpretation of the three Pythagorean centrality measures leads to translate the inequality between them into an inequality between geometric primitives of a hyperrectangle.

\section{Data selection}
\begin{figure}[]
        \centering
        \includegraphics[scale=0.4]{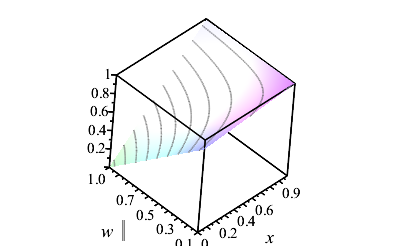}
        \includegraphics[scale=0.4]{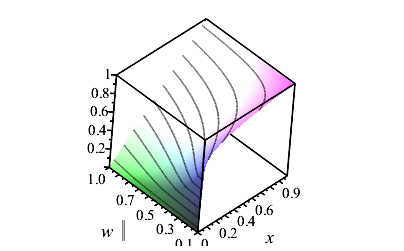}
        \includegraphics[scale=0.4]{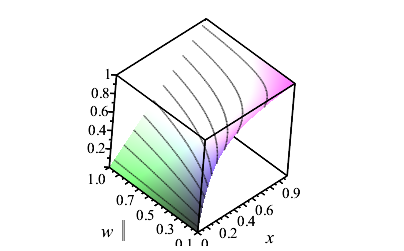}
          \includegraphics[scale=0.4]{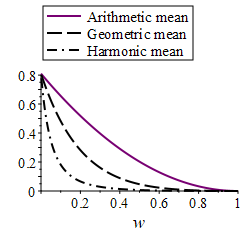}
        \includegraphics[scale=0.4]{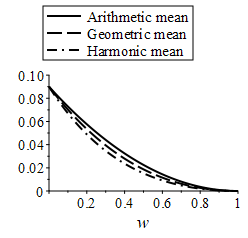}
        \caption{Weighted Pythagorean centrality measures of  1 and $x$ as a function of $x$ and the weight $w$ of $x$. The upper row shows the arithmetic, geometric, and harmonic means respectively. The bottom row shows the velocities of the means when $x=0.1$ (left) and $x=0.7$ (right).}
        \label{fig3M}
    \end{figure}
Remember that the non-negative weight $w_i$ of $x_i$ encodes its relative frequency, the empirical probability of observing it, or its relevance according to a certain criterion. Let us call it, the w-weight. When estimating the mean, the w-weights implement a selection mechanism for measurements, since they indicate how each of them contributes to the estimated mean. Note that the selection of a measurement is not 
taken in the sense of all-or-nothing. Assigning a weight to a measurement is a selection since the measurement's contribution is strengthened or weakened  depending on this weight~\cite{Boutemedjet,Ayech,Bouguila12}. Fig.~\ref{fig3M} displays the three weighted means of two numbers 1 and $x$, by varying $x \in [0,1[$ and the weight $w \in ]0,1]$ of $x$.  For the three weighted means, the more $w$ increases, the more strongly they are attracted by $x$. The difference between the three means is explained by the velocity at which each move towards $x$ as $w$ increases. For a given $w$, the velocity of a mean can be set to the square of the difference between the mean and $x$.  Fig.~\ref{fig3M} presents the velocities of the three means for  $x=0.1$ and $x=0.7$. For a small $x$, the harmonic mean approaches $x$ more quickly than the geometric mean which in turn approaches $x$ more quickly than the arithmetic mean.  This explains why the harmonic mean is the most sensitive to small values. When $x$ increases,  the difference between the velocities of the three means is reduced.

\begin{figure}[]
        \centering
        \includegraphics[scale=0.4]{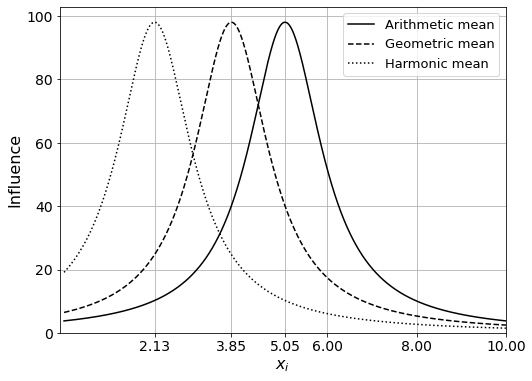}
        \caption{The attraction functions over the three central tendencies obtained by regular sampling of $[0.1, 10]$.}
        \label{fig4}
\end{figure}

\begin{figure}[]
        \centering
        \includegraphics[scale=0.4]{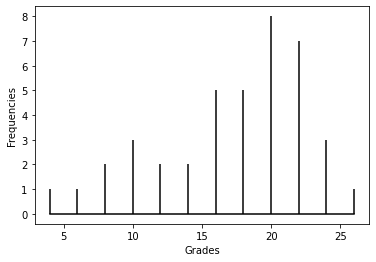}
           \includegraphics[scale=0.4]{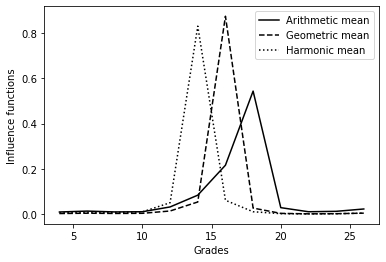}
        \caption{Frequencies of student grades and the three weighted attraction functions over the three central tendencies.}
        \label{Grades}
    \end{figure}

We propose to make explicit another mechanism for selecting measurements  which does not seem to have been made explicit before. The mechanism is expressed in the original space (the measurement space). For the sake of simplicity, first we consider that all measurements have the same w-weight. Let us return to the measurements $x_1 <  \cdots < x_n$. A centrality measure can be seen as a value belonging to the interval $[x_1, x_n]$, where each of the measurements contributes to its identification. The closer a measure is to it, the more influence it has on it to keep it close. This conception implicitly establishes an analogy with the Coulomb's law; i.e. both $x$ and the mean are seen as electric charges whose force of attraction decreases with the square of the distance between them. This analogy  leads to basing attraction  on the Euclidean distance $d(x, \mu)$ between a measurement $x$ and a Pythagorean centrality measure $\mu$ of the measurements.
Since the influence of $x$ on $\mu$ decreases as $d(x, \mu )$ increases, an attraction function can be any decreasing function on $d(x, \mu ) $, such as a Cauchy or Gaussian pdfs. Fig. \ref{fig4} presents an example of the attraction function, which we call the Cauchy attraction function $f_p(x) = \frac{d^2(x_n, x_1)}{d^2(x, \mu_p)+1}$, where $d^2(a,b) = (a-b)^2$ and $p~\in$~\{arithmetic, geometric, harmonic\}. The squared distance $d^2(x_n, x_1)$ is the range of measurements used as the normalization constant. In the denominator, one is added to $d^2()$ in order to avoid  division by zero. 
Fig.~\ref{fig4} depicts the attraction functions obtained by regular sampling of $[0.1, 10]$. The arithmetic, geometric, and harmonic means are   around 5.05,  3.85, and  2.13 respectively. It can be noted that the measurements close to the means have the greatest attraction, regardless of the means used. The attraction of high values is higher in the case of the arithmetic mean. For example, the attraction function of the arithmetic mean at $x=10$ is greater than that of the geometric mean, which in turn is greater than that of the harmonic mean.  Conversely, at $x=0.1$ the attraction function of the harmonic mean is greater than that of the geometric mean which in turn is greater than that of the arithmetic mean.  
Let us examine how the two selection mechanisms (w-weights and attraction function) can be combined. The combination is straightforward because the weighted attraction function can be either  Cauchy   $f_p(x_i) \propto 1/(w_i d^2(x_i, \mu_p)+1)$ or  Gaussian with a non-constant variance $g_p(x_i) \propto exp(-w_i (x_i - \mu_p)^2)$. Fig.~\ref{Grades} shows the student grades and the Cauchy weighted attraction functions for the three means. The observations deduced from Fig.~\ref{fig4} remain valid. Moreover, the data are noisy which explains the difference between  the weighted attraction functions at the maximum points.

The third data selection mechanism is implemented by performing transformations to calculate geometric and harmonic means (see Table~\ref{Pythagore}). Indeed, these two transformations modify the measurements in such a way as to enhance some and weaken the others. For example, the logarithm accentuates the difference between small values and attenuates the difference between large values. Someone may think that the inverse transform neutralizes the effects of the measurements transform. However, the back-transformed arithmetic mean calculated in the transformed space does not cancel the transformation of the measurements because the two transformations are not linear.

To conclude, a weighted mean implements three data selection mechanisms. One of the mechanisms is based on a frequency or any priori knowledge and the two other on the measurement itself. Thus, choosing which Pythagorean centrality measure to use defines how the measurements are considered and which ones are relevant. It follows that the right choice of a mean  can be the cause of appropriate subsequent use of measurements such as accurate predictive analysis and fair decision-making.

\section{Applications}
\subsection{Predictions}
Let's take the example of an airline that needs to estimate the optimal number of additional tickets to sell~\cite{AIRLINES}. These seats are those purchased by passengers, called no-shows, who do not show up for a flight. For estimation purposes, the number of no-shows must be predicted as accurately as possible, because prediction errors lead to financial losses for the airline. A free seat is a loss of earnings and an oversold seat means that a person will be compensated for not boarding. To make the prediction, the airline uses the number of no-shows per flight over the last 200 flights. These data are indicated in  Table~\ref{tab1}. For each  number of no-shows, the table gives the number of flights and the relative frequency of these flights, considered as an empirical probability of the flights. 
\begin{table}[h!]
\centering
 \begin{tabular}{|c|| c| c| c|c|c|c|c|} 
 \hline
 No-shows & 1& 2 & 3 & 4& 5& 6& Total \\ 
 \hline\hline
Number of flights & 70  & 40 & 10 & 20 & 20 & 40 & 200 \\ \hline
 Probability & 0.35 & 0.20 & 0.05 & 0.10 & 0.10 & 0.20 & 1 \\ \hline
 \end{tabular}
   \caption{Number of no-shows on the last 200 flights, the number of flights per number of no-shows, and the empirical probability of the flights. }
 \label{tab1}
\end{table}

%%\begin{table}[h!]
%%\centering
%% \begin{tabular}{|c| c| c|} 
%% \hline
%% Arithmetic mean  & Geometric mean & Harmonic mean  \\ 
%% \hline\hline
%%3.00 &2.34  & 1.83  \\ \hline
%% \end{tabular}
 %%         \caption{Mean values for no-shows}
%%        \label{tab2}
%%\end{table}

%We propose to use each of the Pythagorean centrality measures as a predictor and compare the accuracy of the predictions obtained. The values of these averages are displayed in Table ~\ref{tab2}.
The prediction is made by a consultant whose compensation is linked to the accuracy of the predictions. If $x$ is the expected number of "no-shows" and $n$ the actual number, the compensation is given by the following gain function expressed in a transformed space: 
\begin{equation}
    g(x, n) = 1000 - 30(T(x) - T(n))^2
\end{equation}
As shown in the table~\ref{Pythagore},  $T(a)=a$ gives the arithmetic mean, $T(a)=ln(a)$ the geometric mean, and $T(a)=1/a$ the harmonic mean. Because $n$ is unknown, a way of finding the prediction $x$ is to use the expectation of the gain function, named  the return function.
\begin{equation}
    R(x) =  1000 - \sum\limits^6_{n = 1} 30(T(x) - T(n))^2  pr(n)
  \label{RT}
\end{equation}
%Some values of this function are given in Table~\ref{tab3}.
%\begin{table}[h!]
%\centering
% \begin{tabular}{|c|| c| c| c|c|c|c|c|c|} 
% \hline
%x & 1.00 & 1.83 & 2.00 & 2.34 & 3.00 & 4.00 & 5.00 & 6.00\\ 
% \hline\hline
%R(x) & 763.00  & 841.93 & 853.00 & 869.92 & {\bf 883.00} & 853.00 & 763.00 & 613.00 \\ \hline
% \end{tabular}
%   \caption{Return depending on the predicted number of extra tickets to sell}
%        \label{tab3}
%\end{table}
%The prediction $x=3$ yields the  highest compensation for the consultant. 
Because $T()$ is invertible, maximizing $R(x)$ w.r.t $x$ allows to find the best predictor: 
\begin{equation}
x^* = T^{-1}(\sum\limits^6_{n = 1} T(n)\text{pr}(n))  
\end{equation}
%The critical point $x=AM$ is a maximum because at this point the second derivative  is negative.  In addition to the original data space $x$, let us now express the return function in the two transformed spaces  $ln(x)$ and $1/x$. Let us recall that arithmetic, geometric, and harmonic means  minimize the variance expressed in these three spaces respectively (See Table~\ref{Pythagore}).
Table~\ref{tab4} provides the three return functions, their maximum values, and the best predictors.  For the consultant, the harmonic mean is the best  scenario, since it allows her or him to earn  more. In the absence of ground truth and other economic data, it is difficult to choose the best predictor for the airline. However, the airline could object to the use of a return function and therefore a central tendency measure.
%Sometimes this type of conflict has brought the partners to court to confirm or invalidate the use of a particular mean~\cite{Matthews21}. 

\begin{table}[h!]
\centering
 \begin{tabular}{|c|| c| c| c|c|c|c|c|} 
 \hline
 Return function  $R(x)$ & Best predictor $x^*$ & Best return $R(x^*)$\\ 
 \hline\hline
$1000-30 \sum_{n=1}^6 (x-n)^2 pr(n)$  &AM = 3.00 & 883.00\\ \hline
$1000-30 \sum_{n=1}^6 (ln(x)-ln(n))^2 pr(n) $  &GM=2.34  & 984.26 \\ \hline
$1000-30 \sum_{n=1}^6 (1/x-1/n)^2 pr(n) $  &HM=1.83 & {\bf 996.27}\\ \hline
 \end{tabular}
         \caption{Return functions expressed in the original and transformed spaces of data, their maximum values, and their best predictors}
        \label{tab4}
\end{table}
\subsection{Consumer Price Index}
We briefly describe the case of Canada.
The CPI represents price change by comparing, over time, the cost of the same basket of goods and services. For decades, it was  calculated  using different formulas, including the arithmetic mean, but since 1995, for many products, the geometric mean has been used~\cite{STATCAN: CHAPTER 6}.   The arithmetic mean is still used for three out of 691 products (rents, tuition, and passenger vehicle insurance premiums) on the grounds that they are unlikely to have outliers. It is also used to group products into broader categories, such as food, housing and transport and to calculate the all-items CPI. 
The idea of using the geometric mean instead of the arithmetic mean to estimate CPI was put forward by economist W.~S.~Jevons in 1863~\cite{Jevons}. Later, however, he proposed the use of the harmonic mean instead. The reason for this choice was never explicitly stated in his published work~\cite{Coggeshall}. Let's examine the effect of Pythagorean centrality on the CPI using Canadian data acquired during the two years 2002 and 2017. Note that, the Canadian government tracks a wide variety of products and services. Based on consumer habits, weights are assigned to products indicating their importance to consumers. Table~\ref{tab5} summarizes the weights of the individual product groups and their CPIs in 2017, taking 2002 as the reference period.   Table ~\ref{tab6} shows the overall CPI as a function of Pythagorean centrality measures. As expected, the highest CPI is given by the arithmetic mean and the lowest by the harmonic mean. The latter is 1.3\% lower than the CPI based on the arithmetic mean.
It seems that there are not yet formal arguments to choose the central tendencies to use for the calculation of the CPI. This could open the door to arbitrariness and arguments of a political nature. Nevertheless, a good understanding of data selection could guide the choice of the appropriate central tendencies.
\begin{table}[h!]
\centering
 \begin{tabular}{|c|| c| c|} 
 \hline
 Product group  & Weight & CPI  \\ 
 \hline\hline
Food & 0.1648  & 141.5 \\ \hline
Shelter & 0.2736  & 137.8 \\ \hline
Household operations, & 0.1280  & 121.4 \\ 
furnishings, and equipment &   &  \\ \hline
Clothing and footwear  & 0.0517  & 91.1 \\ \hline
Transportation  & 0.1995  & 133.0 \\ \hline
Health and personal care & 0.0479  & 123.4 \\ \hline
Recreation, education,  and reading & 0.1024  & 111.3 \\  \hline
Alcohol, beverages, tobacco, & 0.0321  & 158.7 \\ 
products, and  cannabis & & \\  \hline
 \end{tabular}
         \caption{Weights and CPI associated with each basket category as of January 2017~\cite{STATCAN: BASKET WEIGHTS,STATCAN: CPI CATEGORY}.}
        \label{tab5}
\end{table}

\begin{table}[h!]
\centering
 \begin{tabular}{|c| c| c|} 
 \hline
 Arithmetic mean  & Geometric mean & Harmonic mean  \\ 
 \hline\hline
130.20 &129.40  & 128.50 \\ \hline
 \end{tabular}
          \caption{All-items CPI given by the three Pythagorean means.}
        \label{tab6}
\end{table}    

\subsection{Ellipse fitting}
\begin{figure}[]
        \centering
        \includegraphics[scale=0.5]{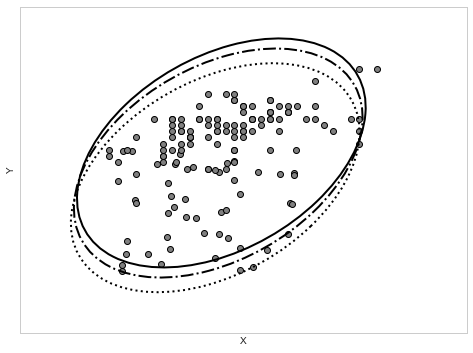}
        \caption{Ellipses obtained by using the arithmetic mean (solid), geometric mean (-.), and harmonic mean (..). }
        \label{Ellipse}
\end{figure}

Fitting ellipses to given points is a problem that arises in  computer graphics, remote sensing, geology, environment,  and statistics, among many others. Different ideas have been implemented, including  least-squares error minimizing and probability density estimation~\cite{Wong12,Chen11,Goumeidane22}. For 2D scatter points,  we want to estimate an enclosing ellipse. An ellipse is defined by its centre and  its principal axes. The centre can be estimated by using a Pythagorean centrality measure of the coordinates of the points. The directions of the principal axes are the eigenvectors of the covariance matrix estimated in the original data space or in the transformed space and transported back to the original space. The lengthening of the ellipse on each of the two axes is determined by the eigenvalues of the covariance matrix. Each of the three means is used for the calculation of both the center of mass and the principal axes of the ellipse.  We can see in Fig.~\ref{Ellipse} that the direction of the principal axes of the ellipses is almost the same, but there is no spatial coincidence of their centres. This is expected because the harmonic mean and the associated ellipse are shifted toward the origin of the data space while the arithmetic mean and the associated ellipse are the farthest.

\section*{Conclusion}
The arithmetic, geometric, and harmonic means make it possible to estimate the centrality of measurable phenomena. The arithmetic mean is the most used. The use of geometric and harmonic means reduces the influence of large measurements because they give  greater influence to  small measurements, which leads to the inequality: AM $>$ GM $>$ HM. Initially, this inequality was displayed in the case of the central tendencies  of two measurements. We have provided a generalization that gives a geometric interpretation of these means in the case of several measurements. Through examples, we have shown the importance of selecting the appropriate mean in the sale of airline tickets, the ellipse estimation, and the estimation of the  Consumer Price Index. We have also explained three data selection mechanisms implemented by Pythagorean means. It should be noted that the choice of the mean has been the subject of several lawsuits. Ultimately, this paper is just a brief overview of the concept of Pythagorean centrality measures. Much remains to be done to understand and better exploit the concept of centrality for data selection. For example, it's hard to summarize and fully understand a set of measurements by using a single centrality measure. Selecting  measurements using different  centrality measures could provide a more complete picture of  data.

\section*{Acknowledgement} \vspace{-2mm}
I would like to thank A. Hajji for  his help in the development of this work. 
%\newpage

\end{document}